\documentclass[12pt]{amsart}
\usepackage{amsmath}
\usepackage{amssymb}
\usepackage{amsfonts}
\usepackage{amsthm}
\theoremstyle{plain}
\newtheorem{theorem}{Theorem}[section]

\newtheorem{proposition}[theorem]{Proposition}

\newtheorem{lemma}[theorem]{Lemma}

\theoremstyle{definition}
\newtheorem{definition}[theorem]{Definition}
\theoremstyle{remark}

\newcommand{\refE}[1]{(\ref{E:#1})}
\newcommand{\refS}[1]{Section~\ref{S:#1}}
\newcommand{\refSS}[1]{Section~\ref{SS:#1}}
\newcommand{\refT}[1]{Theorem~\ref{T:#1}}

\newcommand{\refP}[1]{Proposition~\ref{P:#1}}

\newcommand{\C}{\ensuremath{\mathbb{C}}}
\newcommand{\N}{\ensuremath{\mathbb{N}}}

\renewcommand{\P}{\ensuremath{\mathbb{P}}}

\newcommand{\nl}{\hfill\newline}
\renewcommand{\i}{{\,\mathrm{i}\,}}

\newcommand{\ghm}[1][m]{\Gamma_{hol}(M,L^{#1})}
\newcommand{\gh}{\Gamma_{hol}(M,L)}

\newcommand{\gul}{\Gamma_{\infty}(M,L)}

\renewcommand{\d}{\partial}
\newcommand{\db}{\overline{\partial}}

\newcommand{\zb}{\overline{z}}

\newcommand{\w}{\ensuremath{\omega}}

\newcommand{\volu}{\mathrm{vol}}
\newcommand{\skp}[2]{{\langle #1,#2\rangle}}

\def\pbar{\overline{\partial}}
\def\la{\lambda}
\def\a{\alpha}

\def\Cim{C^{\infty}(M)}
\def\ghm{\Gamma_{hol}(M,L^{m})}

\def\gh{\Gamma_{hol}(M,L)}

\def\gul{\Gamma_{\infty}(M,L)}

\def\Lp{{\mathrm{L}}^2(M,L)}
\def\Lpm{{\mathrm{L}}^2(M,L^m)}

\def\Lqv{{\mathrm{L}}^2(Q,\mu)}

\def\Tfm{T_f^{(m)}}
\def\Tgm{T_g^{(m)}}
\def\Tfgm{T_{\{f,g\}}^{(m)}}
\def\Tma#1{T_{#1}^{(m)}}
\def\Hm{\mathcal {H}^{(m)}}
\def\Hc{{\mathcal {H}}}
\def\d{\partial}
\def\db{\overline{\partial}}

\def\End{\mathrm{End}}

\def\Tr{\operatorname{Tr}}
\def\volu{\operatorname{vol}}
\def\hh{\hat h}
\begin{document}
\vspace*{-1cm}
\hspace*{\fill} Mannheimer Manuskripte 247


\hspace*{\fill} math.QA/9910137

\vspace*{2cm}
\title[Deformation quantization]
{Deformation quantization of compact K\"ahler manifolds
by Berezin-Toeplitz quantization}
\author{Martin Schlichenmaier}
\address{Department of Mathematics and 
  Computer Science, University of Mannheim, D7, 27 \\
         D-68131 Mannheim \\
         Germany}
\email{schlichenmaier@math.uni-mannheim.de}
\begin{abstract}
For arbitrary compact quantizable K\"ahler manifolds it is
shown how a natural formal deformation quantization (star product)
can be obtained via     Berezin-Toeplitz operators.
Results on their semi-classical behaviour (their asymptotic
expansion) 
due to Bordemann, Meinrenken and Schlichenmaier
are used in an essential manner. 
It is shown that the star product is null on constants and fulfills 
parity. A trace is constructed and the relation to
deformation quantization by geometric quantization is given.
\end{abstract}
\date{22.10.99}
\maketitle
{\it \hfill dedicated to the memory of Moshe Flato\qquad}

\section{Introduction}\label{S:intro}
By Bayen, Flato, Fronsdal, Lichnerowicz and 
Sternheimer in 1977
the important concept of quantization given 
by deforming the algebra of functions in ``direction''
of the Poisson bracket was introduced  \cite{BFFLS}.
Clearly the intuitive concept of 
$\hbar$-depending ``deformation'' of classical mechanics
into quantum mechanics was around earlier (e.g. Weyl quantization).
But in their work a mathematically very precise meaning was given 
to it.

Since this time the existence of a
deformation quantization for every symplectic mani\-fold
was established in different ways.
Some of the persons involved were
De Wilde and Lecomte \cite{DeWiLe},
Fedosov \cite{Fed},  and
Omori, Maeda and Yoshioka \cite{OMY}.
Quite recently this was extended to
every Poisson manifold by
Kontsevich \cite{Kont}.
Classification results are also available
\cite{BeCaGu},\cite{Delstar},\cite{NeTS},\cite{Fedbook},\cite{WeXu}.

Even if there is now a very general existence theorem it
is still of importance to study 
deformation quantizations for such manifolds which carry  
additional 
geometric structures.
{}From the whole set of  deformation quantization one is 
looking for one which  keeps
the additional structure.
In this spirit the article deals with the deformation quantization of 
compact quantizable K\"ahler manifolds.
It was shown 1993 
by Bordemann, Meinrenken and Schlichenmaier \cite{BMS} that 
for compact quantizable K\"ahler manifolds  the 
Berezin-Toeplitz quantization has  the correct semi-classical
behaviour (see \refT{approx} below).
Shortly after \cite{BMS} was submitted we had also the result
that by the techniques developed there it was possible to 
construct  a deformation quantization
\cite{DeleMe}.
Details were written up in German \cite{Schlhab} and the result
(with few steps of the proof) appeared in
\cite{SchlBia95}, \cite{SchlGos}.
The complete  proof was not published in English.

Compact K\"ahler manifolds   appear  as phase spaces 
of constrained systems and as 
reduced phase spaces under a  group action.
More recently, they play a rather prominent role in 
Chern-Simons theory, topological and 2-dimensional conformal field theory.
Here typically, 
the phase-spaces to be quantized are 
moduli spaces of certain geometric objects. As examples the 
compactified moduli spaces of
stable holomorphic vector bundles (maybe with additional structures)
on a Riemann surface show up.
The quantum Hilbert spaces  appearing  in this context
are the Verlinde spaces.
 
Encouraged by the recent interest in deformation quantization 
evolving in these fields I found it worthwhile 
to publish the above mentioned results also in English and add some pieces to
it.
By the construction of the deformation quantization 
direct relations to the Berezin-Toeplitz quantization, 
the geometric quantization (via Tuynman's relation),
and asymptotic operator representations are given.
Hence what is presented here is more than just another existence proof.

In the proof the theory of generalized Toeplitz operators
developed by Boutet de Monvel and Guillemin 
\cite{BGTo},\cite{GuCT} is used in an essential manner.
I should mention that in the meantime Guillemin himself
published a sketch of a proof \cite{Guisp} how a deformation 
quantization can be obtained from their theory  
of symbols of Toeplitz operators.

Only for certain special examples of compact K\"ahler manifolds direct
constructions have been known earlier; see results by
Berezin \cite{Bere},
Moreno and Ortega-Navarro \cite{Mor1},\cite{MoOr},
 and Cahen, Gutt and Rawnsley
\cite{CGR2}. 
Recently, 
for all K\"ahler manifolds
(including the noncompact ones)  
the existence of a  
deformation quantization with ``separation of variables''
was shown by Karabegov \cite{Karasep}. 
Separation of variables says essentially that the 
deformation quantization ``respects'' the complex structure.
A classification of all such 
deformation quantizations for a fixed K\"ahler manifold was also given
by Karabegov.
Note that his existence proof is 
on the level of the formal deformation quantization. It does not yield
Hilbert spaces and quantum operators  like in  our 
approach (which in contrast is restricted to the case of quantizable
compact K\"ahler manifolds).
Independently, a similar existence theorem was  proven  by
Bordemann and Waldmann \cite{BW}
along Fedosov's original approach.
Yet another construction was given recently by 
Reshetikhin and Takhtajan 
\cite{ResTak}.

Finally, let me stress the fact, that the 
very essential basics of this work
go back to  joint work 
with Martin Bordemann and Eckhard Meinrenken.
Details have been added by me later on.

The article is organized as follows.
In \refS{set} the geometric set-up is given and the 
main result of this article, the theorem on the construction
of the deformation quantization (\refT{star}) is formulated.
The approximation results from \cite{BMS} are recalled.
In \refS{toe} the  necessary details about the 
Toeplitz structure introduced by Boutet de Monvel and Guillemin
are given.
They are employed in  \refS{proof}
for  the construction of the deformation quantization
(the star product),
i.e. the proof of \refT{star}.
In the concluding 
 \refS{prop} additional properties of the 
star product are discussed.
It is shown that we have $1\star g=g\star 1=g$, i.e. that the 
star product is ``null on constants'' and that it fulfills
the parity condition.
A trace is constructed. 
By a result of Tuynman for compact K\"ahler manifolds the geometric 
quantization can be expressed in terms of the Berezin-Toeplitz quantization.
Using our theorem we see that the geometric quantization 
yields also a star product.
This star product is equivalent to the constructed one. 
The Berezin-Toeplitz star product will be a local star product given by
bidifferential operators. It will have the   
property of ``separation of variables''.
The locality  is not shown in this article. 
The details of the proof are not completely written up
 \cite{KarSchl}.
To be on the safe side one might choose to call 
locality to be conjecturally true.
\section{The Set-Up and the Main Result}\label{S:set}
%
Let $\ (M,\w)\ $ be a compact (complex) K\"ahler manifold.
It should be considered as  phase space manifold $M$
with symplectic form given by the K\"ahler form $\w$.
Denote by $\Cim$ the algebra of (arbitrary often) 
differentiable functions.
Using the K\"ahler form 
one  assigns to every  $f\in\Cim$ its Hamiltonian vector field $X_f$ and 
to every pair of functions $f$ and $g$ the 
Poisson bracket:
\begin{equation}
\label{E:Poi}
\w(X_f,\cdot)=df(\cdot),\qquad  
\{\,f,g\,\}:=\w(X_f,X_g)\ .
\end{equation}
With the Poisson bracket 
$\Cim$ 
becomes a Poisson algebra.

Assume $(M,\w)$ to be quantizable. This says that there exists an
associated quantum line bundle $(L,h,\nabla)$
with holomorphic line bundle $L$,  Hermitian metric  $h$ on $L$ and 
connection $\nabla$ compatible with the metric $h$ and the complex structure
such that the 
curvature  of the line bundle and the K\"ahler form $\w$ of the manifold
are related  as 
\begin{equation}\label{E:quant}
curv_{L,\nabla}(X,Y):=\nabla_X\nabla_Y-\nabla_Y\nabla_X-\nabla_{[X,Y]}
=
-\i\w(X,Y)\ .
\end{equation}
Equation \refE{quant} is called the {\it quantization 
condition}.
If the metric is represented as 
a function $\hh$ 
with respect to 
local complex coordinates  and a 
local holomorphic frame of the bundle
the quantization condition reads as 
$\ \i\db\d\log\hh=\w\ $.

The quantization condition  implies that $L$ is a positive 
line bundle. By the Kodaira embedding theorem   $L$ 
is  ample, which says that  a certain
tensor power $L^{m_0}$ of $L$ is very ample, i.e. 
the global holomorphic sections of  $L^{m_0}$ 
can be used to embed the phase space manifold $M$ 
into projective
space.
Note that the embedding is an embedding as complex manifolds
not as K\"ahler manifolds.
The embedding 
dimension is given by the Hirzebruch-Riemann-Roch formula.
Hence, quantizable compact K\"ahler manifolds are 
as complex manifolds projective
algebraic manifolds.
The converse is also true, see \cite{SchlBia95},\cite{BerSchlcse}.
In the following we will assume $L$ to be  very ample.
If $L$ is not very ample we choose $m_0\in\N$  such that 
the bundle $L^{m_0}$ is very ample and take  
this bundle as quantum line bundle and $m_0\w$  as K\"ahler form 
for $M$. The   underlying complex manifold structure
is not changed.
Please note that for the examples of moduli spaces
mentioned in the introduction  there 
is often  a natural ample or very ample quantum line bundle.

We take the Liouville measure 
$\ \Omega=\frac 1{n!}\w^n\ $ as volume form  on $M$.
On the space of $C^\infty$-sections $\gul$  we have the scalar product
and norm
\begin{equation}
\label{E:skp}
\langle\varphi,\psi\rangle:=\int_M h (\varphi,\psi)\;\Omega\  ,
\qquad
||\varphi||:=\sqrt{\langle \varphi,\varphi\rangle}\ .
\end{equation}
Let $\Lp$ be the  L${}^2$-completion of the space of $C^\infty$-sections
of the bundle $L$   and
$\gh$ be its (due to compactness of $M$) 
finite-dimensional closed subspace of holomorphic
sections.
Let $\ \Pi:\Lp\to\gh\ $ be the projection.
\begin{definition}
For $f\in\Cim$ the {\it Toeplitz operator}  $T_f$
is defined to be
\begin{equation}
T_f:=\Pi\, (f\cdot):\quad\gh\to\gh\ .
\end{equation}
\end{definition}
In words: One takes a holomorphic section $s$ 
and multiplies it with the
differentiable function $f$. 
The resulting section $f\cdot s$ will only be  differentiable.
To obtain a holomorphic section one has to project  it back
on the subspace of holomorphic sections.

The linear map
\begin{equation}
T:\Cim\to \End\big(\gh\big),\qquad  f\to T_f\ ,
\end{equation}
is the  {\it Berezin-Toeplitz quantization map}. 
Because in general
\begin{equation}
T_f\, T_g=\Pi\,(f\cdot)\,\Pi\,(g\cdot)\,\Pi\ne
\Pi\,(fg\cdot)\,\Pi =T_{fg}\ ,
\end{equation}
it is
neither a Lie algebra homomorphism nor
an associative algebra homomorphism.
{}From the point of view of Berezin's approach \cite{Bere}
the operator $T_f$ has as a  contravariant symbol $f$
(see also \cite{SchlBia98} for  relations to 
Berezin's covariant symbols).

This defines a map
from the commutative algebra of functions to a noncommutative
finite-dimensional (matrix) algebra.
The finite-dimensionality is 
due to compactness of $M$.
A lot of classical information will get lost. To recover this
information one should consider not just the bundle $(L,\nabla,h)$ alone but
all its tensor powers $(L^m,\nabla^{(m)},h^{(m)})$
and apply the above constructions for every $m$.
Note that if
 $\hat h$ corresponds to the metric $h$ w.r.t. a holomorphic frame $s$
of the bundle $L$ then $\hat h^m$ corresponds to the metric $h^{(m)}$ 
w.r.t. to the frame $s^{\otimes m}$ for the bundle $L^m$.
In this way one obtains a family of
matrix algebras and a family of maps
\begin{equation}
\Tma {}:\Cim\to \End\big(\ghm\big),\qquad  f\to \Tma f\ .
\end{equation}
This infinite family should in some sense ``approximate'' the
algebra $\Cim$.(See \cite{BHSS}  and the discussion on
strict quantization below.)
for a definition  of such an approximation.)
Indeed this family has the correct semi-classical behaviour 
as is expressed in \refT{approx} below.

It also allows to construct a deformation quantization.
A deformation quantization is given by a {\it star product}.
I will use both terms interchangeable.
To fix the notation and the factors of $\i$ 
let me recall the definition of a star product.
Let $\mathcal{A}=\Cim[[\nu]]$ be the algebra of formal power
 series in the
variable $\nu$ over the algebra $\Cim$. A product $\star$
 on $\mathcal {A}$ is 
called a (formal) star product if it is an
associative $\C[[\nu]]$-linear product such that
\begin{enumerate}
\item
\qquad $\mathcal{ A}/\nu\mathcal {A}\cong\Cim$, i.e.\quad $f\star g
 \bmod \nu=f\cdot g$,
\item
\qquad $\dfrac 1\nu(f\star g-g\star f)\bmod \nu = -\i \{f,g\}$,
\end{enumerate}
where $f,g\in\Cim$. 
We can also write 
\begin{equation}
\label{E:cif}
 f\star g=\sum\limits_{j=0}^\infty C_j(f,g)\nu^j\ ,
\end{equation}
with
$ C_j(f,g)\in\Cim$. The $C_j$ should be  $\C$-bilinear in $f$ and $g$.
The conditions 1. and 2.  can 
be reformulated as 
\begin{equation}
\label{E:cifa}
C_0(f,g)=f\cdot g,\qquad\text{and}\qquad
C_1(f,g)-C_1(g,f)=-\i \{f,g\}\ .
\end{equation}
The aim of this article
is to show the following 
\begin{theorem}
\label{T:star}
There exists a unique (formal) star product on $\Cim$
\begin{equation}
f * g:=\sum_{j=0}^\infty \nu^j C_j(f,g),\quad C_j(f,g)\in
C^\infty(M),
\end{equation}
in such a way that for  $f,g\in\Cim$ and for every $N\in\N$  we have
with suitable constants $K_N(f,g)$ for all $m$
\begin{equation}
\label{E:sass}
||T_{f}^{(m)}T_{g}^{(m)}-\sum_{0\le j<N}\left(\frac 1m\right)^j
T_{C_j(f,g)}^{(m)}||=K_N(f,g) \left(\frac 1m\right)^N\ .
\end{equation}
\end{theorem}
This theorem has been proven immediately after 
\cite{BMS} was finished. It has been announced in \cite{SchlBia95},
\cite{SchlGos}
and the proof was written up in German in  \cite{Schlhab}.
In \refS{proof} I will supply the proof.

Instead of writing \refE{sass} we will sometimes use the more
intuitive notation 
\begin{equation}
\label{E:expans}
T_{f}^{(m)}\cdot T_{g}^{(m)}
\quad \sim\quad \sum_{j=0}^\infty\left(\frac 1m\right)^j
T_{C_j(f,g)}^{(m)}
\qquad (m\to\infty)\ .
\end{equation}
The asymptotics should always be understood in the above 
precise sense.

In the proof the  results 
expressed in the following theorem are needed.
Denote by $||f||_\infty$ the sup-norm of $\ f\ $ on $M$ and by 
$||\Tfm||=\sup_{s\in\ghm, s\ne 0}\frac {||\Tfm s||}{||s||}$
 the operator norm on $\ghm$.
\newpage
\begin{theorem}
\label{T:approx}
[Bordemann, Meinrenken, Schlichenmaier]

(a) For every  $\ f\in \Cim\ $ there exists $C>0$ such that   
\begin{equation}
||f||_\infty +\frac Cm
\le||\Tfm||\le ||f||_\infty\ .
\end{equation}
In particular, $\lim_{m\to\infty}||\Tfm||= ||f||_\infty$.

(b) For every  $f,g\in \Cim\ $ 
\begin{equation}
\label{E:dirac}
||m\i[\Tfm,\Tgm]-\Tfgm||\quad=\quad O(\frac 1m)\quad
\mathrm{as}\quad m\to\infty
\ .
\end{equation}

(c) For every  $f,g\in \Cim\ $ 
\begin{equation}
||\Tfm\Tgm-T^{(m)}_{f\cdot g}||\quad=\quad O(\frac 1m)\quad
\text{as}\quad m\to\infty
\ .
\end{equation}
\end{theorem}
These results are contained in  Theorem  4.1, 4.2,
resp. in Section 5 in \cite{BMS}.
Note that part (c) also follows from \refE{sass} for $N=1$ and
generalizes trivially to finitely many functions.

Our result does not prove a strict deformation quantization in the 
sense of Rieffel \cite{Rief}.
But it is a {\it strict quantization} (see 
for the definition \cite{LandTop},\cite{Riefque}).
Let $I:=\{\frac 1m\mid m\in\N\}\cup \{0\}$ be the topological space  with
topology coming from the real line.
It has $0$ as accumulation point.
To every $\hbar\in I$, $\hbar\ne 0$,
i.e $\hbar=1/m$, one assigns
the algebra
$A_{1/m}:=\End(\ghm)$ with $||.||_{1/m}$ the operator norm and
 to $0$ one assigns the algebra  $A_0:=\Cim$ with norm $||.||_0=|.|_{\infty}$.
The map $\hbar\to T^{(1/\hbar)}_f$, with  $T^{(\infty)}_f:=f$ defines by
\refT{approx} a continuous field of $C^*$-algebras on the
family 
$(A_\hbar)_{\hbar\in I}$.
{}From \refE{Tconj} follows that $T$ respects conjugation.
By \refE{dirac} the additional condition for a strict
quantization is also fulfilled.
Due to the compactness of $M$ the maps 
$T^{(1/\hbar)}_f$ for $\hbar\ne 0$ 
are never injective. Hence the strict quantization 
is not faithful at a fixed level $\hbar$, only in the limit 
$\hbar\to\infty$.

In \cite{BHSS} and \cite{BMS} the notion of $L_{\alpha}$, resp.
$gl(N)$, resp. $su(N)$ quasi-limit was used for this concept.
It was conjectured in \cite{BHSS} that for every 
compact K\"ahler manifold the Poisson algebra of function 
is a $gl(N)$ quasi-limit. This was proven in \cite{BMS}.
This result is of special interest in the theory of membranes.

There is another geometric concept of quantization, the {\it geometric
quantization} introduced by Kostant and Souriau.
But for compact K\"ahler manifolds 
due to Tuynman \cite{TuyQ} (see also \cite{BHSS} for a coordinate independent 
proof) 
they have the same
semi-classical behaviour
\begin{equation}
\label{E:tuyn}
Q_f^{(m)}=\i\cdot T_{f-\frac 1{2m}\Delta f}^{(m)}\ . 
\end{equation}
Here 
 $Q_f^{(m)}$ is the well-known  operator
of
geometric quantization 
(with respect 
to the quantum line bundle $L^m$)  corresponding to the
prequantum operator
$\ P_f^{(m)}= -\nabla_{X_f^{(m)}}^{(m)}+\i f\cdot id\ $ and 
K\"ahler polarization.
K\"ahler polarization means 
$Q_f^{(m)}=\Pi^{(m)}P_f^{(m)}\Pi^{(m)}$
with the projectors 
\begin{equation}
\label{E:prom}
\Pi^{(m)}:\Lpm\to\ghm\ .
\end{equation}
In \refE{tuyn}
$\Delta$ is the Laplacian with respect to the K\"ahler metric given by
$\omega$.
In \refS{prop} I will show that this allows to define a deformation
quantization via the operators of geometric quantization. It will be 
equivalent to the Berezin-Toeplitz deformation quantization.
\section{Toeplitz Structure}\label{S:toe}
%
In \cite{BMS} the set-up for the proof of the 
approximation results was given. Here I use the same 
setting.
Let me recall for further reference the main definitions.
A  more detailed exposition can be found in 
\cite{Schlhab}.
Take $\ (U,k):=(L^*,h^{-1})\ $ the dual of the quantum line bundle,
$Q$ the unit circle bundle inside $U$ (with respect to the metric $k$)  and
$\tau: Q\to M$ the projection.
Note that for the projective space 
with quantum line bundle the hyperplane section bundle  $H$,
the bundle $U$ is just the tautological
bundle. Its  fibre over the point $z\in\P^N(\C)$ consists of
the line in $\C^{N+1}$ which is represented by $z$. In particular,
for the projective space 
the total space of $U$ with  the zero section removed can be identified
with $\C^{N+1}\setminus\{0\}$.
The same picture remains true for the via the very ample quantum
line bundle in projective space embedded manifold $M$.
The quantum line bundle will be the pull-back of $H$ 
(i.e. its restriction to the embedded manifold) and its
dual is the pull-back of the tautological bundle.

In the following we use $E\setminus 0$ to denote the total space of
the vector bundle $E$ with the image of the zero section removed. 
Starting from the function
$\hat k(\la):=k(\la,\la)$ on $U$ we define
 $\tilde a:=\frac {1}{2\i}(\d-\pbar)\log \hat k$ on
$U\setminus 0$ (with respect to the complex
structure on $U$) and denote by $\alpha$
its restriction   to $Q$.
Now $d\a=\tau^*\w$ (with $d=d_Q$) and $\mu=\frac 1{2\pi}\tau^*\Omega
\wedge \a$ is a volume form
on $Q$. With respect to this form we take the L${}^2$-completion $\Lqv$
of the space of functions on $Q$.
The generalized Hardy space $\Hc$ is the closure of the
functions in  $\Lqv$ which can be extended to
holomorphic functions on the whole
disc bundle.
The generalized Szeg\"o projector is the projection
\begin{equation}
\label{E:szproj}
\Pi:\Lqv\to \Hc\ .
\end{equation}
 By the natural circle action
$Q$ is a $S^1$-bundle and the tensor powers of $U$ can be
viewed as associated bundles. The space $\Hc$ is preserved
by this action.
It can be decomposed into eigenspaces 
$\Hc=\prod_{m=0}^\infty \Hm$ where
 $c\in S^1$ acts on $\Hm$ as multiplication
by $c^m$.
Sections of $L^m=U^{-m}$ can be identified with functions $\phi$ on $Q$ which
satisfy the equivariance condition
$\phi(c\la)=c^m\phi(\la)$.
It turns out that  this 
identification is an isometry.
Recall that $\Lpm$ has a scalar product given in an corresponding
way to \refE{skp}. 
Restricted to the holomorphic objects
we obtain an isometry
\begin{equation}
\ghm\cong\Hm\ .
\end{equation}
There is the notion of Toeplitz structure
$(\Pi,\Sigma)$ as developed by Boutet de Monvel  and
Guillemin
in \cite{BGTo},\cite{GuCT}.
What is needed from there are only the following facts.
$\Pi $ is the  Szeg\"o projector  \refE{szproj}.
The second object is  the submanifold 
\begin{equation}
\Sigma=\{\;t\alpha(\lambda)\;|\;\lambda\in Q,\,t>0\ \}\ \subset\  T^*Q
\setminus 0
\end{equation}
of the tangent bundle of $Q$ 
defined with the help of the 1-form $\alpha$.
They showed that it is a 
symplectic submanifold.
A (generalized) Toeplitz operator of order $k$ is  an operator
$A:\Hc\to\Hc$ of the form
$\ A=\Pi\cdot R\cdot \Pi\ $ where $R$ is a
pseudodifferential operator
($\Psi$DO) of order $k$ on
$Q$.
The Toeplitz operators build a ring.
The  symbol of $A$ is the restriction of the
principal symbol of $R$ (which lives on $T^*Q$) to $\Sigma$.
Note that $R$ is not fixed by $A$, but
Guillemin and Boutet de Monvel showed that the  symbols
are well-defined and that they obey the same rules as the
symbols of   $\Psi$DOs.
In particular we have the following relations
\begin{equation}
\label{E:symbol}
\sigma(A_1A_2)=\sigma(A_1)\sigma(A_2),\qquad
\sigma([A_1,A_2])=\i\{\sigma(A_1),\sigma(A_2)\}_\Sigma.
\end{equation}
In our context  only  two Toeplitz operators appear:
\nl
(1) The generator of the circle action
gives the  operator $D_\varphi=\dfrac 1{\i}\dfrac {\partial}
{\partial\varphi}$. It is an operator of order 1 with symbol $t$.
It operates on $\Hm$ as multiplication by $m$.
\nl
(2) For $f\in\Cim$ let $M_f$ be the  operator on
$\Lqv$ 
corresponding to multiplication with $\tau^*f$.
We set 
\footnote{
There should be no confusion with the operator
$T_f=T_f^{(1)}$ introduced above.}
$\ T_f=\Pi\cdot M_f\cdot\Pi:\Hc\to\Hc\ $.
Because $M_f$ is constant along the fibres of $\tau$, $T_f$ 
commutes with the circle action.
Hence
$\ T_f=\prod\limits_{m=0}^\infty\Tfm\ $,
where $\Tfm$ denotes the restriction of $T_f$ to $\Hm$.
After the identification of $\Hm$ with $\ghm$ we see that these $\Tfm$
are exactly the Toeplitz operators  $\Tfm$ introduced in \refS{set}.
In this sense  $T_f$ is called  the global Toeplitz operator and
the $\Tfm$ the local Toeplitz operators.
$T_f$ is an operator of order $0$.
Let us denote by
$\ \tau_\Sigma:\Sigma\subseteq T^*Q\to Q\to M$ the composition
then we obtain for the symbol  $\sigma(T_f)=\tau^*_\Sigma(f)$.
\section{Proof of \refT{star}}\label{S:proof}
Let the notation be as  in the last section.
In particular, let $T_f$ be the 
Toeplitz operator, $D_\varphi$ the operator of rotation, and 
$\Tfm$, resp. $(m\cdot)$ their projections on the eigenspaces 
$\Hm\cong\ghm$.

\noindent
\subsection{The definition of the  $C_j(f,g)\in \Cim$}
The construction is done inductively in such a way
that 
\begin{equation}
\label{E:Asum}
A_N=D_\varphi^N T_fT_g- \sum_{j=0}^{N-1}
D_\varphi^{N-j}T_{C_j(f,g)}
\end{equation}
is always a  Toeplitz operator of order zero.
The operator $A_N$ is $S^1$-invariant, i.e.
$D_\varphi\cdot A_N=A_N\cdot D_\varphi$. 
Because it is of order zero his symbol 
is a function on $Q$. By the $S^1$-invariance the symbol is even 
given by (the pull-back of) a function on $M$.
Denote this function to be the next element $C_N(f,g)$
in the star product.
By construction the 
operator  $\ A_N-T_{C_N(f,g)}\ $ is of order   $-1$ and 
$A_{N+1}=D_\varphi(A_N-T_{C_N(f,g)})$ is of order  0 
and hence exactly of the form given 
in \refE{Asum}.
The  induction starts with
\begin{gather}
A_0=T_fT_g, \qquad\mathrm{and}\
\\
\sigma(A_0)=\sigma(T_f)\sigma(T_g)=
\tau^*_{\Sigma}(f)\cdot 
\tau^*_{\Sigma}(g)
=
\tau^*_{\Sigma}(f\cdot g)\ .
\end{gather}
Hence, $C_0(f,g)=f\cdot g$ as required.
\newline
It remains to show statement \refE{sass}
about the asymptotics.
As an operator of order zero on a compact manifold
$A_N$ is bounded. Hence the same is true for all its restrictions
$A_N^{(m)}$ to $\Hm$.
If we calculate them we obtain
\begin{equation}
\label{E:azwi}
||m^N\Tfm\Tgm-\sum_{j=0}^{N-1}m^{N-j}T^{(m)}_{C_j(f,g)}||
=||A_N^{(m)}||\le ||A_N||\ .
\end{equation}
After dividing by $m^N$ Equation \refE{sass} follows.
Bilinearity is clear.

\noindent
\subsection{The Poisson structure}
The relation  $C_0(f,g)=f\cdot g$ was proven above.
To show the 2. formula in 
\refE{cifa} we write explicitly  \refE{azwi} 
for  $N=2$ and the pair of functions
$(f,g)$:
\begin{equation}
||m^2 \Tfm\Tgm-m^2 T^{(m)}_{f\cdot g}-m T^{(m)}_{C_1(f,g)}||
\le K\ .
\end{equation}
A corresponding expression is obtained for the pair  $(g,f)$.
If we subtract both operators inside of the norm 
we obtain (with the triangle inequality and suitable $K'$)
\begin{equation}
||m^2 (\Tfm\Tgm-\Tgm\Tfm)-m (T^{(m)}_{C_1(f,g)}-T^{(m)}_{C_1(g,f)})||
\le K'\ .
\end{equation}
Dividing by $m$ and multiplying with $i$ we obtain
\begin{equation}
||m\i [\Tfm,\Tgm]-T^{(m)}_{\i\big(C_1(f,g)-C_1(g,f)\big)}||
=O(\frac 1m)\ .
\end{equation}
Using the asymptotics given by \refT{approx}(b) for 
the commutator we get
\begin{equation}
|| T^{(m)}_{\{f,g\}-\i\big(C_1(f,g)-C_1(g,f)\big)}||=O(\frac 1m)\ .
\end{equation}
Taking the limit for $m\to\infty$ and using 
\refT{approx}(a) we get
\begin{equation}\ ||\{f,g\}-\i(C_1\big(f,g)-C_1(g,f)\big)||_\infty=0\ . 
\end{equation}
Hence $\ \{f,g\}=\i(C_1(f,g)-C_1(g,f))$.
This shows  \refE{cifa}.

\noindent
\subsection{The uniqueness}
It is proven by induction using the asymptotics
\refE{sass}.
Let $C_j(f,g)$ and $\tilde C_j(f,g)$
be two such systems of bilinear maps fulfilling the 
required properties.
Assume  $C_j=\tilde C_j$ for 
$j\le N-2$. If we subtract the corresponding expressions
in \refE{sass} and use the fact that $T^{(m)}$ is linear
we obtain
\begin{equation}
||\frac 1{m^{N-1}}T^{(m)}_{(C_{N-1}(f,g)-\tilde C_{N-1}(f,g))}||
\le \frac K{m^N}\ .
\end{equation}
Hence,
\begin{equation}
\lim_{m\to\infty}
||T^{(m)}_{(C_{N-1}(f,g)-\tilde C_{N-1}(f,g))}||
=0\ .
\end{equation}
With  \refT{approx} (a) it follows  
$\ C_{N-1}(f,g)=\tilde C_{N-1}(f,g)$.
The induction starts with  $N=1$.
But here  $C_0(f,g)=\tilde C_0(f,g)=f\cdot g$
is required.

\noindent
\subsection{The associativity}
The proof employs the associativity of the operators 
used to construct the star product  and again 
\refT{approx} (a).
The relation 
$\ f\star(g\star h)=(f\star g)\star h\ $ 
can be rewritten in relations for the maps $C_j$:
\begin{equation}
\label{E:sasso}
\sum_{l=0}^k C_l(f,C_{k-l}(g,h))=
\sum_{l=0}^k C_l(C_{k-l}(f,g),h)\ .
\end{equation}
{}From  \refT{approx} (a) we know
$
f=g\quad\longleftrightarrow\quad
\lim_{m\to\infty}||T_f^{(m)}-T_g^{(m)}||
=0$.
Hence it is enough to apply the  Toeplitz operator $T^{(m)}$ 
to the relation \refE{sasso} and study the asymptotics of 
$\ T_{\mathrm{left\ hand\ side}}- T_{\mathrm{right\ hand\ side }}\ $.
\nl
This is done by induction over $k$.
\nl
$k=0:$ $C_0(f,C_0(g,h))=C_0(C_0(f,g),h)$ is true because
$C_0(f,g)=f\cdot g$.
\nl
Assume the claim to be true up to  level
 $k-1$.
The equation \refE{sass} for   $0\le r\le k$ multiplied by 
 $m^r$ ($N=r+1$) yields
\begin{equation}
\label{E:zwi2}
T_{C_r(f,g)}^{(m)}=m^rT_f^{(m)}T_g^{(m)}-
\sum_{s=0}^{r-1}m^{r-s}T_{C_s(f,g)}
+O(\frac 1m)\ .
\end{equation}
Here the symbol $O(\frac 1m)$ 
is shorthand for the statement that the difference of the 
operators on the left and on the right is an operator 
whose norm behaves like
$O(\frac 1m)$ 
for $m\to\infty$.
In particular we obtain for  $l=0,1,\ldots,k$
\begin{equation}
T_{C_l(f,C_{k-l}(g,h))}^{(m)}=m^lT_f^{(m)}T_{C_{k-l}(g,h)}^{(m)}-
\sum_{s=0}^{l-1}m^{l-s}T_{C_s(f,C_{k-l}(g,h))}
+O(\frac 1m)\ .
\end{equation}
Summation over $l$ yields 
\begin{equation}
T^{(m)}_{\mathrm{l.h.s.}}=
\sum_{l=0}^k m^l
T_f^{(m)}T_{C_{k-l}(g,h)}^{(m)}
-\sum_{l=0}^k
\sum_{s=0}^{l-1}m^{l-s}T^{(m)}_{C_s(f,C_{k-l}(g,h))}
+O(\frac 1m)\ .
\end{equation}
The second sum can be rewritten as 
\begin{equation}
-\sum_{r=1}^k m^r\sum_{l=r}^k T^{(m)}_{C_{l-r}(f,C_{k-l}(g,h))}
=-
\sum_{r=1}^k m^rT^{(m)}_{\sum_{s=0}^{k-r}C_s(f,C_{k-r-s}(g,h))}\ .
\end{equation}
For such  sums we know by induction that \refE{sasso} 
is valid. The same is done for the right hand side.
If we subtract 
$T^{(m)}_{\mathrm{r.h.s.}}$ from
$T^{(m)}_{\mathrm{l.h.s.}}$, 
it remains
\begin{equation}
\sum_{l=0}^k m^l T_f^{(m)} T_{C_{k-l}(g,h)}^{(m)}
-
\sum_{l=0}^k m^l T_{C_{k-l}(f,g)}^{(m)}T_h^{(m)}
+O(\frac 1m)\ .
\end{equation}
By splitting the first sum into 
$l=0$ and $l\ge 1$ and using for the
$l=0$ term  the asymptotic
\refE{zwi2} we obtain
\begin{gather}
m^0\left(
T_f^{(m)}m^k\Tgm T_h^{(m)}-\sum_{s=0}^{k-1} m^{k-s}T_f^{(m)} 
T_{C_s(g,h)}^{(m)}+O(\frac 1m)\right)
+\sum_{l=1}^k m^l
 T_f^{(m)} 
T_{C_{k-l}(g,h)}^{(m)}
\\
=
m^k T_f^{(m)}(\Tgm T_h^{(m)})+O(\frac 1m)\ .
\end{gather}
A corresponding expression follows for the second sum.
As difference remains
\begin{equation}
m^k (T_f^{(m)}(\Tgm T_h^{(m)})-(T_f^{(m)}\Tgm) T_h^{(m)})
+O(\frac 1m)\ .
\end{equation}
Now we ended up with operators which are clearly associative,
The operator coming with 
the $m^k$ term vanishes.
Hence associativity follows from 
\refT{approx} (a).
\qed
\section{Additional Properties}\label{S:prop}
The introduced star product has important  properties.
\subsection{Unit}
The unit of the algebra $\Cim$, the constant function 1, will 
also be the unit in the star product.
Such star products are sometimes called to have the
property "null on constants".
\begin{proposition}
For the above introduced star product we have
\begin{equation}
\label{unit}
1\star g= g\star 1=g\ .
\end{equation}
Equivalently,
\begin{equation}
\label{unitc}
C_k(1,g)=C_k(g,1)=0,\qquad \mathrm{for}\quad k\ge 1\ .
\end{equation}
\end{proposition}
\begin{proof}
For $f\equiv 1$ we have  $\ T_f\equiv id\ $, resp. $\ \Tfm\equiv id\ $.
Also  $C_0(1,g)=g=C_0(g,1)$. Further with  \refE{Asum} 
\begin{equation}
A_1=D_\varphi T_f T_g- D_\varphi T_{fg}=
D_\varphi T_g-D_\varphi T_g\ .
\end{equation}
Hence the symbol of $A_1$ vanishes. But this implies
$C_1(1,g)=0=C_1(g,1)$.
The claim follows by trivial  induction from \refE{Asum} .
\end{proof}
\subsection{Parity}
A star product is said to fulfill the parity condition if
\begin{equation}
\label{E:pari}
\overline{f\star g}=\overline{g}\star \overline{f}\ .
\end{equation}
Considering the formal parameter to be real ($\overline{\nu}=\nu$)
this is  equivalent to  
\begin{equation}
\label{E:paric}
\overline{C_k(f,g)}={C_k(\overline{g},\overline{f})},\quad  k\ge 0\  .
\end{equation}
We will show
\begin{proposition}
\label{P:pari}
The above introduced star product fulfills parity.
\end{proposition}
\begin{lemma}
\begin{equation}
\label{E:Tconj}
{T_f^{(m)}}^*=T_{\overline{f}}^{(m)}\ .
\end{equation}
\end{lemma}
\begin{proof}
Take any $s,t\in \ghm$.
For the scalar product we calculate
($\Pi^{(m)}$ is the projector defined  in \refE{prom})
\begin{equation}
\skp {s}{T_f^{(m)}t}=
\skp {s}{\Pi^{(m)}ft}=
\skp {s}{ft}=
\skp {\overline{f}s}{t}=
\skp {T_{\overline{f}}^{(m)}s}{t}\ .
\end{equation}
Hence the claim.
\end{proof}
\begin{proof} (\refP{pari})
Recall that the identification of the sections
of $L^m$ with equivariant functions on the circle bundle $Q$
is an isomorphy.
Hence  the definition of  adjoint operators agree.
For the global Toeplitz operator we obtain
$T_f^*=\prod_{m=0}^\infty T_{\overline{f}}^{(m)}=T_{\overline{f}}$.
The star product 
\begin{equation}
\overline{g}\star \overline{f}
=\sum_{j=0}^\infty \nu^jC_j(\overline{g},\overline{f})
\end{equation}
is given 
via the asymptotic expansion of 
\begin{equation}
T_{\overline{g}}^{(m)}\cdot
T_{\overline{f}}^{(m)}
=
{T_{g}^{(m)}}^*\cdot 
{T_{f}^{(m)}}^*
=
({T_{f}^{(m)}}\cdot 
{T_{g}^{(m)}})^* \ .
\end{equation}
For the asymptotic expansion 
of the last expression we have 
\begin{equation}
({T_{f}^{(m)}}\cdot 
{T_{g}^{(m)}})^*
\quad\sim\quad
\sum_{j=0}^\infty \left(\frac {1}{m}\right)^j T_{C_j(f,g)}^{(m)*}=
\sum_{j=0}^\infty \left(\frac {1}{m}\right)^j 
T_{\overline{C_j(f,g)}}^{(m)}\ .
\end{equation}
But this is the complex conjugate of the asymptotic expansion 
which defines $f\star g$.
This shows \refE{pari}
\end{proof}
By the parity condition we have on $\Cim[[\nu]]$ an
anti-involution given by pointwise complex conjugation on the functions,
and by considering the
formal parameter to be real ($\overline{\nu}={\nu}$). 
\subsection{Locality and Separation of Variables}
Recall that a star product is local if for all $f,g\in\Cim$ the support
 $supp\ C_j(f,g)$ is contained in $supp\ f\cap supp\ g$ for all
$j\in\N_0$.
Using Peetre's theorem and the fact 
that the $C_j$ are bilinear 
this implies that for a local star product
the $C_j$ can be given by bidifferential operators.

Using the fact that the projection operators $\Pi^{(m)}$ can be 
expressed 
with the help of Berezin-Rawnsley's coherent states and the fact that
the coherent states are ``localizing'' for $m\to\infty$,  
it is possible to show that the above star product is local.
The details are not written up 
\cite{KarSchl}. Hence I prefer to call this statement 
conjecturally true.

{}From
the expected behaviour  it follows   that 
the star product can be restricted to 
open subsets and defines  compatible star products there.
For such star products Karabegov introduced 
the notion of star products with {\it separation of variables} \cite{Karasep}
(Bordemann and Waldmannn \cite{BW} called them star products
of Wick type).
In our convention this reads as 
$f\star k=f\cdot k$ and $k\star g=k\cdot g$ for (locally defined)
holomorphic functions $g$, antiholomorphic functions $f$
and arbitrary functions $k$.
The above introduced star product will 
be a star product with separation of variables
(assuming locality).

More precisely, we expect that the $C_j$ are bidifferential operators
of degree $(j,j)$ with only holomorphic derivatives in the
first entry and only antiholomorphic derivatives in the
second entry.
See \refSS{example} for examples.
\subsection{Trace}
\begin{proposition} [Bordemann, Meinrenken, Schlichenmaier]
Let $f\in\Cim$ and let $n=\dim_{\C}M$. Denote the 
trace on $\End(\ghm)$ by $\Tr^{(m)}$  then
\begin{equation}
\Tr^{(m)}\,(T^{(m)}_{f})
        =m^n\left(\frac 1{\volu (\P^{n}(\C))}
\int_M f\, \Omega +O(m^{-1})\right)\ .
\end{equation}
\end{proposition}
This result can also be found in \cite{BMS}.
There it was given only with a hint of its proof.
Because it is central for the following let me 
give the details.
\begin{proof}
Let us start with a real valued  $f$.
Then the operator $T_f$ and the components  $T_f^{(m)}$
are self adjoint (see Equation \refE{Tconj}).
Let $d(m)=\dim\Hm$    and let 
$\lambda_1^{(m)},\lambda_2^{(m)},\ldots,\lambda_{d(m)}^{(m)}$ 
be the eigenvalues of the restriction  of 
$T_f$ on $\Hm$.
In particular, these are also the eigenvalues 
of $\Tfm$ on
$\ghm$.
Following \cite{BGTo}  ($n=\dim_{\C} M$)
let  
\begin{equation}
\mu_m=\frac 1{m^n}\sum_{i=1}^{d(m)}\delta(\lambda-\lambda_i^{(m)})
\end{equation}
be the discrete spectral measure.
By Theorem~13.13 of \cite{BGTo}
it converges weakly to the limit
measure 
\begin{equation}
\mu(g)=\gamma_M\int_M g(f(z))\,\Omega(z)
\end{equation}
with a universal constant $\gamma_M$ only
depending on the manifold $M$.
An important intermediate result there  is
the asymptotic expansion (Equation 13.13 in \cite{BGTo})
\begin{equation}
\label{E:specas}
\mu_m(g)\quad\sim\quad\sum_{r=-n}^\infty a_r(g)m^{r+n}\ .
\end{equation}
For $g\equiv 1$ we obtain
\begin{equation}
\label{E:spectr}
\frac 1{m^n}\sum_{i=1}^{d(m)}\lambda_i^{(m)}
=\frac 1{m^n}\Tr^{(m)} \Tfm=
\gamma_M\int_M f\,\Omega+O(\frac 1m)\ .
\end{equation}
To calculate $\gamma_M$ we evaluate \refE{spectr} 
for  $f\equiv 1$ (i.e. $\Tfm=id$) 
and obtain 
\begin{equation}
\gamma_M=\frac {\dim\ghm}{m^n\cdot\volu(M)}+O(\frac 1m)\ .
\end{equation}
Note that  (see  p.113 and Thm 5.22 in \cite{MumCP})
\begin{equation}
\dim\ghm=\frac {\volu(M)}{\volu(\P^{n}(\C))}\cdot m^n+O(m^{n-1})\ .
\end{equation}
Hence 
\begin{equation}
\gamma_M={\volu(\P^{n}(\C))}^{-1}.
\end{equation}
In particular the coefficient depends only on the dimension of
$M$.
This shows the claim for real valued $f$. For complex valued  
$f$ it follows from linearity by considering real and imaginary 
part separately.
In \cite{MumCP} for $M$ the restriction
$\omega_{FS|M}$ of the Fubini-Study K\"ahler form was used
to define the volume. Here we have to work with the form $\omega$.
Because the deRham classes of both forms coincide and because 
the K\"ahler forms are closed the volume will be the same.
\end{proof}

{}From \refE{specas} follows the 
asymptotic expansion for $m\to\infty$ (see also \cite{BPUspec})
\begin{equation}
\label{E:tras}
\Tr^{(m)}(T_f^{(m)})\quad\sim\quad
m^n\left(\sum_{j=0}^\infty \left(\frac 1{m}\right)^j\tau_j(f)\right),
\quad \mathrm{with}\quad \tau_j(f)\in\C\ .
\end{equation}
We define the $\C[[\nu]]$-linear map
\begin{equation}
\label{E:trace}
\Tr:\Cim[[\nu]]\to \nu^{-n}\C[[\nu]],\quad
\Tr f:=\nu^{-n}\sum_{j=0}^\infty\nu^j\tau_j(f),
\end{equation}
such that for $f\in\Cim$ the $\tau_j(f)$ are given
by the asymptotic expansion \refE{tras}
and for arbitrary elements by $\C[[\nu]]$-linear extension.
\begin{proposition}
The map $\Tr$ is a trace, i.e. we have
\begin{equation}
\label{E:trsym}
\Tr (f\star g)=\Tr (g\star f)\ .
\end{equation}
\end{proposition}
\begin{proof}
By $\C[[\nu]]$-linearity it is enough to show this for $f,g\in\Cim$.
The element  $f\star g-g\star f$ is given by
the asymptotic expansion
 of $T_f^{(m)}\cdot T_g^{(m)}-T_g^{(m)}\cdot T_f^{(m)}$.
Hence $\Tr(f\star g-g\star f)$ is given by the expansion of 
\begin{equation}
\Tr^{(m)}(T_f^{(m)}\cdot T_g^{(m)}-T_g^{(m)}\cdot T_f^{(m)})\ .
\end{equation}
But for every $m$ this vanishes.
Hence \refE{trsym} follows.
\end{proof}

\subsection{Examples}
\label{SS:example}
For the sphere $S^2$, resp. $\P^1(\C)$ with K\"ahler form
\begin{equation}
\omega=\frac {\i}{(1+z\zb)^2}\; dz\wedge d\zb\ ,
\end{equation}
and the hyperplane bundle  as  quantum line-bundle 
explicit calculations 
\footnote{
Not following the lines of the proof in \refS{proof} 
but working with a basis of the
sections of the bundles.
}
 of the author (not published)
yield (using  $T^{(m)}(f):=T_f^{(m)}$)
\begin{equation}
\label{E:goa}
\lim_{m\to\infty}||m\left(
T^{(m)}(f) 
T^{(m)}(g)-
T^{(m)}(fg)\right)+
T^{(m)}\bigg((1+z\zb)^2\frac {\partial f}{\partial z}
\frac {\partial g}{\partial \zb}\bigg)||=0\ .
\end{equation}
This implies 
\begin{equation}
\label{E:zwi3}
C_1(f,g)=-(1+z\zb)^2\frac {\partial f}{\partial z}
\frac{\partial g}{\partial \zb}\ .
\end{equation}
For the case of Riemann surfaces of genus  $g\ge 2$
more than  half of the article \cite{KlLeb} 
by Klimek and Lesniewski 
deals with the proof of the  fact corresponding
to \refE{goa}.
In the realization of the Riemann surface $M$ as quotient space
$\{z\in\C\mid |z|<1\}/G$ with $G$ a Fuchsian subgroup of 
$\mathrm{SU}(1,1)$ acting by fractional linear 
transformations one takes as K\"ahler form 
the $\mathrm{SU}(1,1)$ invariant form
\begin{equation}
\omega=\frac {2\i}{(1-z\zb)^2}dz\wedge d\zb\ .
\end{equation}
The corresponding quantum line bundle is the canonical
line bundle, i.e. the bundle whose local sections are the
holomorphic differentials.
{}From their results  follows
\begin{equation}
C_1(f,g)=-\frac 12 (1-z\zb)^2\frac {\partial f}{\partial z}
\frac {\partial g}{\partial \zb}\ .
\end{equation}
\subsection{Deformation quantization via geometric quantization}
Via Tuynman's relation \refE{tuyn} the operator $Q^{(m)}$ of geometric
quantization corresponding to the function $f$ can be expressed by the
Toeplitz operator $T^{(m)}$ corresponding to the function
$f-\frac 1{2m}\Delta f$.
\refT{star} shows that the asymptotic expansion
\begin{equation}
\label{E:geoas}
Q^{(m)}_f\cdot Q^{(m)}_g
\quad\sim\quad
\sum_{j=0}^\infty\left(\frac 1m\right)^j
Q^{(m)}_{D_j(f,g)}\ ,
\end{equation}
with suitable $D_j(f,g)\in\Cim$
is well-defined in the precise sense as expressed in the theorem there.
We set
\begin{equation}
f\star_G g:=\sum_{j=0}^\infty \nu^jD_j(f,g)\ .
\end{equation}
The first two terms calculate as
\begin{gather}
D_0(f,g)=f\cdot g 
\\
D_1(f,g)=C_1(f,g)+\frac 12\left(
\Delta(f\cdot g)-\Delta f\cdot g-f\cdot \Delta g\right),
\end{gather}
where the $C_j$ are the coefficients of the Berezin-Toeplitz star product.
In particular the conditions \refE{cifa} are fulfilled for the $D$'s.
Hence, this defines indeed a star product.

In fact more is valid. If we introduce the linear maps
\begin{equation}
B^{(m)}(f):=f-\frac 1{2m}\Delta f,
\end{equation}
and the $\C[[\nu]]$-linear map induced by
\begin{equation}
B(f):=f-\nu\frac {\Delta}{2} f=
(id-\nu\frac {\Delta}{2})f 
\end{equation}
on $\Cim[[\nu]]$ we can rewrite \refE{geoas}
\begin{equation}
T^{(m)}_{B^{(m)}(f)}\cdot T^{(m)}_{B^{(m)}(g)}
\quad\sim\quad
\sum_{j=0}^\infty\left(\frac 1m\right)^j
T^{(m)}_{B^{(m)}(D_j(f,g))}\ .
\end{equation}
Taking the asymptotics we get
\begin{equation}
\label{equiv}
B(f)\star B(g)=B(f\star_G g) \ .
\end{equation}
Note that $B(f) \mod \nu=f$, $B(1)=1$ and that
$B$ is invertible. The inverse is given  by
\begin{equation}
B^{-1}=id+\sum_{k=1}^\infty \frac 1{2^k}\nu^k\Delta^{k}\ .
\end{equation}
Recall that two star products (over the same manifold) are 
equivalent if there exists a $\C[[\nu]]$-algebra isomorphism inducing
the identity on the zero order part.
This implies
\begin{proposition}
The star product of geometric quantization is equivalent to the star product
of Berezin-Toeplitz quantization.
\end{proposition}

\bigskip
\subsection*{Acknowledgements}
As the basics of the presented work goes back to  joint work
with M. Bordemann and E. Meinrenken 
it is a pleasure for me to thank them for all their inspirations.
I also like to thank the Erwin Schr\"odinger 
International Institute for Mathematical Physics 
in Vienna for its hospitality,
and F. Haslinger, P. Michor and H. Upmeier for the 
invitation to participate in one of
the activities at the institute.
During the stay there   the main part of this article was written up.

\begin{thebibliography}{10}

\bibitem{BFFLS}
Bayen,~F., Flato,~ M., Fronsdal,~C., Lichnerowicz,~A., and Sternheimer,~D.:
  {Deformation theory and quantization, Part {I}}, 
   \textit{Lett. Math. Phys.}
  \textbf{1} (1977), 521--530.
 {Deformation theory and quantization, Part {II} and {III}}, \textit{Ann.
  Phys.} \textbf{111} (1978), 61--110, 111--151.

\bibitem{BerSchlcse}
  Berceanu,~St., and Schlichenmaier,~M.: {Coherent state embeddings, polar
  divisors and {C}auchy formulas}, \textit{math.QA/9902066,
 Mannheimer Manuskripte Nr.   232.}

\bibitem{Bere}
  Berezin,~F.A.: {Quantization}, \textit{Math. USSR-Izv.} \textbf{8} (1974),
  1109--1165.
  {Quantization in complex symmetric spaces}, \textit{Math. USSR-Izv.}
  \textbf{9} (1975), 341--379.

\bibitem{BeCaGu}
  Bertelson,~M., Cahen,~M., and Gutt,~S.: {Equivalence of star products},
  \textit{Classical Quantum Gravity} \textbf{14} (1997), A93--A107.

\bibitem{BHSS}
 Bordemann,~M., Hoppe,~J., Schaller,~P., and Schlichenmaier,~M.:
   {$gl(\infty)$
  and geometric quantization}, \textit{Commun. Math. Phys}.
  \textbf{138} (1991),
  209--244.

\bibitem{BMS}
  Bordemann,~M., Meinrenken,~E., and Schlichenmaier,~M.:{{T}oeplitz
  quantization of {K\"{a}}hler manifolds and $gl(n), n\to\infty$ limits},
  \textit{Commun. Math. Phys.} \textbf{165} (1994), 281--296.

\bibitem{BW}
   Bordemann,~M., and Waldmann,~St.:{A {F}edosov star product of the {W}ick
  type for {K}{\"a}hler manifolds}, \textit{Lett. Math. Phys.} \textbf{41} (1997),
  243--253.

\bibitem{BPUspec}
  Borthwick,~D., Paul,~D., and Uribe,~A.:{Semiclassical spectral estimates for
  {T}oeplitz operators}, \textit{ Ann. Inst. Fourier, Grenoble}
   \textbf{48} (1998),
  1189--1229.

\bibitem{BGTo}
 Boutet~de Monvel,~L., and Guillemin,~V.:\emph{The spectral theory of {T}oeplitz
  operators. {A}nn. {M}ath. {S}tudies, {N}r.99}, Princeton University Press,
  Princeton, 1981.

\bibitem{CGR2}
  Cahen,~M., Gutt,~S., and Rawnsley,~J.: \emph{Quantization of {K\"{a}}hler
  manifolds {II}}, Trans. Amer. Math. Soc. \textbf{337} (1993), 73--98.

\bibitem{DeWiLe}
De~Wilde,~M., and Lecomte,~P.B.A.:{Existence of star products and of formal
  deformations of the {P}oisson-{L}ie algebra of arbitrary symplectic
  manifolds}, \textit{Lett. Math. Phys.} \textbf{7} (1983), 487--496.

\bibitem{DeleMe}
   Deligne,~P.: {Letter to {B}ordemann, {M}einrenken and {S}chlichenmaier
  ({M}arch 1994). Answer by {M}einrenken}, 1994.

\bibitem{Delstar}
\bysame, {Deformation de l'algebre des fonctions d'une variete
  symplectique: Comparaison entre {F}edosov et {D}e {W}ilde,{L}ecomte}, 
  \textit{Sel.   Math., New Ser.} \textbf{1} (1995), 667--697.

\bibitem{Fedbook}
Fedosov,~B.V.: \emph{Deformation quantization and index theory}, Akademie
  Verlag, Berlin, 1996.

\bibitem{Fed}
  Fedosov,~B.V.: {Deformation quantization and asymptotic operator
  representation}, \textit{Funktional Anal. i. Prilozhen.}
   \textbf{25} (1990), 184--194.
   {A simple geometric construction of deformation quantization}, 
   \textit{ J.  Diff. Geo.} \textbf{40} (1994), 213--238.

\bibitem{GuCT}
   Guillemin,~V: {Some classical theorems in spectral theory revisited},
   in \textit{
  Seminars on singularities of solutions of linear partial differential
  equations}, Ann. Math. Studies, Nr.91 (L.~H{\"o}rmander, ed.), Princeton
  University Press, 1979, pp.~219--259.

\bibitem{Guisp}
\bysame, {Star products on pre-quantizable symplectic manifolds}, 
  \textit{Lett.
  Math. Phys.} \textbf{35} (1995), 85--89.


\bibitem{Karasep}
  Karabegov,~A.: {Deformation quantization with separation of variables on
  a {K}{\"a}hler manifold}, \textit{Commun. Math. Phys.}
   \textbf{180} (1996), 745--755.



\bibitem{KarSchl}
  Karabegov,~A., and Schlichenmaier,~M.: \emph{in preparation}.

\bibitem{KlLeb}
  Klimek,~S., and Lesniewski,~A., {Quantum {R}iemann surfaces: {II}. {T}he
  discrete series}, \textit{Lett. Math. Phys.} \textbf{24} (1992), 125--139.

\bibitem{Kont}
  Kontsevich,~M.: {Deformation quantization of {P}oisson manifolds,{I}},
  q-alg/9709040.

\bibitem{LandTop}
 Landsman,~N.P.: \emph{Mathematical topics between classical and quantum
  mechanics}, Springer, Berlin, Heidelberg, New York, 1998.

\bibitem{Mor1}
  Moreno,~C.: {$*$-products on some {K}{\"a}hler manifolds}, 
 \textit{Lett. Math.   Phys.} \textbf{11} (1986), 361--372.

\bibitem{MoOr}
 Moreno,~C., and Ortega-Navarro,~P.:{{$*$}-products on {$D^1(\C)$}, {$S^2$}
  and related spectral analysis}, \textit{Lett. Math. Phys.} \textbf{7} (1983),
  181--193.

\bibitem{MumCP}
 Mumford,~D.: \emph{Complex projective varieties}, Springer, Berlin, Heidelberg,
  New York, 1976.

\bibitem{NeTS}
  Nest,~R., and Tsygan,~B.: {Algebraic index theory}, 
  \textit{Commun. Math. Phys.}
  \textbf{172} (1995), 223--262.
   {Algebraic index theory for families},  \textit{Advances in Math}
   \textbf{113}
  (1995), 151--205.

\bibitem{OMY}
  Omori,~H., Maeda,~Y., and Yoshioka,~A.:{Weyl manifolds and deformation
  quantization}, \textit{Advances in Math} \textbf{85} (1991), 224--255.
  {Existence of closed star-products}, \textit{Lett. Math. Phys.}
  \textbf{26} (1992), 284--294.

\bibitem{ResTak}
 Reshetikhin,~N., and Takhtajan,~L.:{Deformation quantization of {K\"a}hler
  manifolds}, math/9907171.

\bibitem{Rief}
   Rieffel,~M.A.: {Deformation quantization and operator algebras},
   in: \textit{ Operator
  theory/Operator algebras and applications}, Proc. Sympos. Pure Math. 51
  (W.~Arveson and R.~Douglas, eds.), AMS, 1990, pp.~411--423.



\bibitem{Riefque}
 Riefel,~M.: {Questions on quantization}, 
  in: \textit{Operator algebras and operator
  theory. Proceedings of the international conference, Shanghai, China, July
  4--9, 1997}, (Liming et~al., ed.), AMS, 1998, pp.~315--326.


\bibitem{SchlBia95}
  Schlichenmaier,~M.: 
 {{B}erezin-{T}oeplitz quantization of compact {K\"{a}}hler
  manifolds}, in: \textit{
 Quantization, Coherent States and Poisson Structures, Proc.
  XIV'th Workshop on Geometric Methods in Physics (Bia\l owie\.za, Poland, 9-15
  July 1995)} (A.~Strasburger, S.T. Ali, J.-P. Antoine, J.-P. Gazeau, and
  A.~Odzijewicz, eds.), Polish Scientific Publisher PWN, 1998, q-alg/9601016,
  pp.~101--115.


\bibitem{Schlhab}
\bysame, \emph{Zwei {A}nwendungen algebraisch-geometrischer {M}ethoden in der
  theoretischen {P}hysik: {B}erezin-{T}oeplitz-{Q}uantisierung und globale
  {A}lgebren der zweidimensionalen konformen {F}eldtheorie}, Habilitation
  Thesis, 1996.


\bibitem{SchlGos}
\bysame, {Deformation quantization of compact {K\"{a}}hler manifolds via
  {B}erezin-{T}oeplitz operators}, 
   in: \textit{Proceedings of the XXI Int. Coll. on Group
  Theoretical Methods in Physics (15-20 July 1996, Goslar, Germany)}, (H.-D.
  Doebner, P.~Nattermann, and W.~Scherer, eds.), World Scientific, 1997,
  pp.~396--400.


\bibitem{SchlBia98}
\bysame, {{B}erezin-{T}oeplitz quantization and {B}erezin
  symbols for arbitrary compact {K}{\"a}hler manifolds}, to appear in 
  \textit{Rep. on   Math. Phys.}, math.QA/9902066.


\bibitem{TuyQ}
  Tuynman,~G.M.:{Quantization: Towards a comparison between methods}, 
  \textit{Jour.
  Math. Phys.} \textbf{28} (1987), 2829--2840.

\bibitem{WeXu}
 Weinstein,~A., and Xu,~P.: {Hochschild cohomology and characteristic classes
  for star products},in:
  \textit{ Geometry and Differential equations}, (A.~Khovanskij
  et~al., eds.), AMS, 1998, pp.~177--194.

\end{thebibliography}
\providecommand{\bysame}{\leavevmode\hbox to3em{\hrulefill}\thinspace}

\end{document}